\documentclass[a4paper,11pt]{article}
\usepackage{amssymb}
\usepackage{amsmath}
\usepackage{amsfonts}
\usepackage[mathscr]{eucal}
\usepackage{enumerate}
\usepackage{amscd}
\usepackage{indentfirst}
\usepackage{enumerate}
\usepackage{txfonts, textcomp}
\newtheorem{thm}{Theorem}[section]
\newtheorem{prop}[thm]{Proposition}
\newtheorem{dfn}[thm]{Definition}
\newtheorem{lem}[thm]{Lemma}
\newtheorem{ex}{Example}[section]
\numberwithin{equation}{section}

\allowdisplaybreaks[4]
\title{Notes on Geometric Morita equivalence \\of twisted Poisson manifolds}
\author{Yuji HIROTA\thanks{yhirota@rs.noda.tus.ac.jp}
\\Quantum Bio-Informatics Center\\
Research Institute for Science and Technology\\Tokyo University of Science}
\date{}
\begin{document}
\maketitle 

\begin{abstract}
This paper is devoted to the study of Morita equivalence for twisted Poisson manifolds. We review some Morita 
invariants and prove that integrable twisted Poisson manifolds which are gauge equivalent are Morita equivalent. 
Moreover, we introduce the notion of weak Morita equivalence and show that if two twisted Poisson manifolds are 
weak Morita equivalent, there exists a one-to-one correspondence between their twisted symplectic leaves. 
\end{abstract}

\section{Introduction}
Geometric Morita theory is one of the interesting topics in Poisson geometry. 
The geometric notion of Morita equivalence 
was introduced by Xu, P.~(\cite{xu1},\cite{xu2}and \cite{xu3}) 
on the basis of algebraic Morita equivalence. 
Morita equivalence is first introduced by Morita, K. in \cite{morita}. 
He gave a necessary and sufficient condition for representation categories of two rings to be equivalent: 
two rings have equivalent categories of left modules if and only if there exists an equivalence bimodule for rings. 
Ring theoretical Morita equivalence is generalized to the theory of $C^*$-algebras by Rieffel, M. 
\cite{rieffel1},\cite{rieffel2}. 
Morita equivalence of $C^*$-algebra is useful in studying some $C^*$-algebras. Also, Morita equivalent 
$C^*$-algebras share many properties, such as equivalent categories of Hermitian left modules, 
isomorphic $K$-group, and so on. 
$C^*$-algebras are the quantum objects; in contrast, Poisson manifolds are the classical one. 
Morita equivalence for integrable Poisson manifolds 
and (quasi-) symplectic groupoids ware introduced by Xu as the classical analogue of this equivalence relation. 
Geometric Morita equivalence plays an important role in Poisson geometry as Morita equivalence of $C^*$-algebras does. 
There exist some invariants under Morita equivalence such as the representation 
categories of symplectic realizations, fundamental groups and the first Poisson cohomology groups~
(see Ginzburg, V. L. and Lu, J.-H. \cite{lu} and \cite{xu2}). 
And furthermore, the theory of geometric Morita equivalence is related to momentum map theory. 
Since Morita equivalence establishes an equivalence of representation categories, 
we are provided with the notion of equivalence for momentum map theories. 
It is shown that some known correspondence of momentum map theories can be described by Morita equivalences \cite{xu3}.
On the basis of Xu's work, 
the author introduce the notion of Morita equivalence for integrable twisted Poisson manifolds \cite{yuji1}, \cite{yuji2}. 
As for Poisson manifolds, Morita equivalent twisted Poisson manifolds have isomorphic fundamental groups, isomorphic 
first cohomology groups and equivalent categories of modules.  
Morita equivalence is applied only for integrable (twisted) Poisson manifolds. 
To remedy this defect, we will introduce refined version of Morita equivalence and discuss it in this paper.

The paper is organized as follows: In Section 2 we study the basic properties of twisted Poisson manifolds and discuss 
the relation with Lie algebroids. 
Section 3 begin with the review of Morita equivalence discussed in \cite{yuji1} and \cite{yuji2}. 
The latter part of this section 
deals with Dirac structures and a gauge transformation. After that, 
We will prove that gauge equivalence of integrable twisted Poisson manifolds implies Morita equivalence. 
In Section 4, we introduce the notion of weak Morita equivalence of (twisted) Poisson manifolds 
and show that 
Morita equivalence implies weak Morita equivalence. Furthermore, we define a bijective correspondence 
between the twisted symplectic leaves of $P_1$ and those of $P_2$ when twisted Poisson manifolds $P_1$ and $P_2$ 
are weak Morita equivalent. 

Finally, we note that smooth manifolds appeared in this paper are assumed to be connected. We denote by $\varGamma(E)$ 
the set of smooth sections of a vector bundle $E\to M$. 

\section{Preliminaries}
\subsection{Twisted Poisson manifolds}
Twisted Poisson manifolds first appeared in the study of string theory by 
Park. J.-S. \cite{park} and Klim\v{c}\'{\i}k, C. and Strobl, T. \cite{WZW}, and 
treated mathematically by \v{S}evera, P. and Weinstein, A. \cite{severa}. 
We start by recalling the definition of a twisted Poisson manifold.

A twisted Poisson manifold is a smooth manifold $P$ equipped with a bivector field $\Pi$ and a closed $3$-form 
$\phi$ on $P$ which satisfy the following equation:
\begin{equation}\label{eq1}
\frac{1}{2}[\Pi,\,\Pi] \,=\, \wedge^3\Pi^\sharp(\phi),
\end{equation}
where, $[\cdot,\,\cdot]$ means a Schouten-Nijenhuis bracket and $\wedge^3\Pi^\sharp(\phi)$ is a linear map from 
$\varGamma(\wedge^3T^*P)$ to $\varGamma(\wedge^3TP)$ 
induced from the natural homomorphism 
$\Pi^\sharp : T^*M\to TM$ given by 
$\beta\bigl(\Pi^\sharp(\alpha)\bigr) = \langle\Pi,\,\beta\wedge\alpha\rangle$. 
Namely, for any 
$\alpha,\,\beta,\,\gamma \in \varGamma(T^*P)$, $\wedge^3\Pi^\sharp(\phi)$ is defined as 
\[\wedge^3\Pi^\sharp(\phi)(\alpha,\,\beta,\,\gamma) := \phi\bigl(\Pi^\sharp(\alpha),\,\Pi^\sharp(\beta),\,\Pi^\sharp(\gamma)
\bigr). \] 

We call $\Pi$ a twisted Poisson bivector.
\begin{ex}{\rm 
{\rm (Poisson manifolds)}~Let $(P,\,\Pi)$ be a Poisson manifold. For a closed $3$-form $\phi$ on $P$ such that 
$\wedge^3\Pi^\sharp(\phi)=0$, it holds that $[\Pi,\,\Pi] = 0 = \wedge^3\Pi^\sharp(\phi)$. Therefore, $(P,\,\Pi,\,\phi)$ is a 
twisted Poisson manifold. 
}
\end{ex} 

\begin{ex}{\rm 
Let $A$ be the set of elements $\{x_1,\,x_2,\,x_3,\,x_4\}\subset \mathbb{R}^4$ which satisfy $x_1=0$ or $x_3=0$. 
The closed 3-form $\phi = \bigl((1/x_3^2)dx_2 - (1/x_1^2)dx_4\bigr)\wedge dx_1\wedge dx_3$ on $\mathbb{R}^4\setminus A$ 
and the bivector $\Pi = x_3(\partial /\partial x_1)\wedge (\partial/\partial x_2) 
\,+\, x_1(\partial /\partial x_3)\wedge (\partial/\partial x_4)$ satisfy the condition {\rm (\ref{eq1})}. 
In other words, $(\mathbb{R}^4\setminus A,\,\pi,\,\phi)$ is a twisted Poisson manifold. 
}
\end{ex} 

Given a $\phi$-twisted Poisson manifold $(P,\,\Pi)$, one can define a bilinear skew-symmetric map $\{\cdot,\,\cdot\}$ 
on $C^\infty(P)$ and a vector field on $P$ by 
 \[
  \{f,\,g\} := \langle \Pi,\,df\wedge dg\rangle,\qquad 
        H_f := \Pi^\sharp(df),\quad \bigl(\forall f,\,g\in C^\infty(P)\bigr).
 \]
The vector field $H_f$ determined by $f\in C^\infty(P)$ is called the Hamiltonian vector field of $f$. 
It is easy to verify that the map $\{\cdot,\,\cdot\}$ satisfies the Leibniz identity. 
By using the bracket and the Hamiltonian vector fields, the formula (\ref{eq1})~ can be written as

\begin{equation}\label{Jacobi}
 \bigl\{\{f,\,g\},\,h\bigr\} + \bigl\{\{g,\,h\},\,f\bigr\} + \bigl\{\{h,\,f\},\,g\bigr\} + \phi(H_f,\,H_g,\,H_h) = 0.
\end{equation}

Conversely, if a bilinear skew-symmetric map $\{\cdot,\,\cdot\}:C^\infty(P)\times C^\infty(P)\to C^\infty(P)$
and a closed 3-form $\phi\in \varGamma(\wedge^3T^*P)$ 
satisfy the Leibniz identity and 

\begin{equation}\label{anomaly}
 \left\{\{f,g\},h\right\} \,+\, \left\{\{g,h\},f\right\} \,+\, \left\{\{h,f\},g\right\} \,=\, 
                     \bigl\langle \{f,\,\cdot\}\wedge\{g,\,\cdot\}\wedge\{h,\,\cdot\},\,\phi\bigr\rangle, 
\end{equation}

then $\{\cdot,\,\cdot\}$ arises from a 2-vector field $\Pi$ given by 
\[\langle \Pi,\,df\wedge dg\rangle = \{f,\,g\},\quad (\forall f,\,g\in C^\infty(P)). \]
Furthermore, it can be verified that $\Pi$ and $\phi$ satisfy the formula (\ref{eq1}).
In consequence, we can define 
a twisted Poisson manifold as a smooth manifold $P$ together with a closed 3-form $\phi\in \varGamma(\wedge^3T^*P)$ and 
a bilinear skew-symmetric map $\{\cdot,\,\cdot\}:C^\infty(P)\times C^\infty(P)\to C^\infty(P)$ satisfy the equation 
(\ref{anomaly})~ and the Leibniz identity. 

The following proposition can be shown by simple calculation. 
\begin{prop}\label{eq2} For any $f,\,g,\,h\in C^\infty(P)$, it holds that 
 \[\bigl(\,[H_f,\, H_g] + H_{\{f,g\}}\,\bigr)h = \phi(H_f,\,H_g,\,H_h).\]
\end{prop}
({\it Proof}\,)~Using (\ref{Jacobi}), we have 
\begin{align*}
\bigl([H_f,\, H_g] + H_{\{f,g\}}\bigr)h &= H_f(H_gh) - H_g(H_fh) + H_{\{f,g\}}h\\
                                        &= \bigl\{ \{h,\, g\}, f\bigr\} - \bigl\{ \{h,\, f\}, g\bigr\} + 
                                                                                      \bigl\{h, \{f,\, g\}\bigr\}\\ 
                                        &= -\Bigl(\bigl\{ \{g,\, h\}, f\bigr\} + \bigl\{ \{h,\, f\}, g\bigr\} 
                                                                               + \bigl\{ \{f,\, g\}, h\bigr\}\Bigr)\\
                     &= \phi(H_f,\,H_g,\,H_h).\qquad\qquad\qquad\qquad\qquad\qquad\qquad\qquad\qquad\qquad\quad \Box
\end{align*}

\begin{dfn}
For a closed 3-form $\psi$ on a smooth manifold $S$, $\psi$-twisted symplectic form is 
a non-degenerate 2-form $\omega\in \varGamma(\wedge^2T^*S)$ such that $d\omega=\psi$. 
A smooth manifold equipped with a $\psi$-twisted symplectic form is called a $\psi$-twisted symplectic manifold.
\end{dfn}

The non-degeneracy of a $\psi$-twisted symplectic form $\omega$ implies that 
the natural homomorphism $\omega^\flat:TS\to T^*S,\quad X\mapsto i_X\omega$ is an isomorphism, 
where $i_X\omega$ means the contraction of $\omega$ by $X$. 
Therefore, given a smooth function $f$ on $S$, we can define its Hamiltonian vector field $H_f$ by 
$i_{H_f}\omega = df$. Moreover, as for symplectic manifolds, we can define a bracket on 
$S$ as $\{f,\,g\}:=\omega(H_f,\,H_g)$. Then it is verified that the bracket $\{\cdot,\,\cdot\}$ obtained from 
$\omega$ and $\psi$ satisfy the equation (\ref{Jacobi})~and the Leibniz identity. That is, a twisted symplectic 
manifold is a twisted Poisson manifold. 

\vspace{0.5cm}
Let $(P_i,\,\Pi_i,\,\phi_i)~(i=1,2)$ be twisted Poisson manifolds and $J:P_1\to P_2$ a smooth map. 
The smooth map $J$ is called a twisted Poisson map if, for any $x\in P_1$, the following formula holds:
\begin{equation}\label{twPoisson map}
(\Pi_2^\sharp)_{J(x)} \,=\, (dJ)_x\circ \Pi_1^\sharp\circ (dJ)_x^\ast .
\end{equation}
By using the bracket, a twisted Poisson map $J:P_1\to P_2$ can be written as 
\begin{equation}\label{twPoisson map2}
\{f,\,g\}_2\circ J \,=\, \{J^*f,\,J^*g\}_1,\quad \bigl(\forall f,\,g\in C^\infty(P_2)\bigr),
\end{equation}
where $\{\cdot,\,\cdot\}_i~(i=1,2)$ mean the brackets induced from $\Pi_i$. 

\begin{prop}\label{twPmap}
Let $(P_i,\,\Pi_i,\,\phi_i)~(i=1,2)$ be twisted Poisson manifolds. If a smooth map $J:P_1\to P_2$ is a twisted 
Poisson map, then for any $\alpha,\,\beta,\,\gamma\in T_{J(x)}^*P_2~(\forall x\in P_1)$, 
\begin{align*}
 (\phi_1)_{J(x)}\Bigl(\Pi_1^\sharp\bigl(\alpha\!\circ\!(dJ)_x\bigr),\,&\Pi_1^\sharp\bigl(\beta\!\circ\!(dJ)_x\bigr),\,
     \Pi_1^\sharp\bigl(\gamma\!\circ\!(dJ)_x\bigr)\Bigr)\\  
  &= (J^*\phi_2)_x\Bigl(\Pi_1^\sharp\bigl(\alpha\!\circ\!(dJ)_x\bigr),\,\Pi_1^\sharp\bigl(\beta\!\circ\!(dJ)_x\bigr),\,
     \Pi_1^\sharp\bigl(\gamma\!\circ\!(dJ)_x\bigr)\Bigr).
\end{align*}
\end{prop}

This proposition is easily shown by using (\ref{Jacobi}) and (\ref{twPoisson map2}).  
We remark that $J^*\phi_2=\phi_1$ does not hold in general even if $\Pi_1^\sharp$ is an isomorphism. 

\begin{dfn}
A twisted symplectic realization {\rm (}t.s.realization for short{\rm )} of a twisted Poisson manifold $P$ 
is a twisted symplectic manifold $S$ together with a twisted Poisson map $J:S\to P$. 
\end{dfn}

Analogously, we define an anti-twisted symplectic realization (anti-t.s.realization, for short) of a twisted 
Poisson manifold $(P,\,\Pi,\,\phi)$ as a 
twisted symplectic manifold $S$ together with a twisted Poisson map $J':S\to \overline{P}$, where $\overline{P}$ 
means a twisted Poisson manifold $(P,\,-\Pi,\,-\phi)$. 
\subsection{Lie algebroids of twisted Poisson manifolds}
If $(P,\,\Pi,\,\phi)$ is a twisted Poisson manifold, then the cotangent bundle $T^*P\to P$ carries a 
Lie algebroid structure whose anchor map is the natural anchor map 
$\Pi^\sharp:T^*P\to TP,\,\beta(\Pi^\sharp\alpha) = \langle\Pi,\,\beta\wedge\alpha\rangle$ 
and whose Lie bracket is 
\begin{equation}\label{bracket}
[\alpha,\,\beta]_\phi := \, \mathcal{L}_{\sharp \alpha}\beta \,-\, \mathcal{L}_{\sharp \beta}\alpha \,+\, 
    d\bigl(\Pi (\alpha,\,\beta)\bigr) \,-\, \phi(\sharp\alpha,\,\sharp\beta,\,\cdot\,)
\end{equation}
where we denote by $\mathcal{L}_X\omega$ the Lie derivative of $\omega$ by $X$. 

Let us verify that the Jacobi identity. 
For any $f,\,g,\,h\in C^\infty(P)$, 
\begin{align}\label{algebroid1}
0 &= d\Big(\bigl\{\{f,\,g\},\,h\bigr\} + \bigl\{\{g,\,h\},\,f\bigr\} 
                   + \bigl\{\{h,\,f\},\,g\bigr\} + \phi(H_f,\,H_g,\,\cdot\,)\Bigr)\notag \\ 
  &= d\bigl\{\{f,\,g\},\,h\bigr\} + d\{\{g,\,h\},\,f\bigr\} 
                   + d\bigl\{\{h,\,f\},\,g\bigr\} + d\bigl(\phi(H_f,\,H_g,\,\cdot\,)\bigr)\notag \\ 
  &= \bigl[[df,\,dg],\,dh\bigr] + \bigl[\phi(H_f,\,H_g,\,\cdot\,),\,dh\bigr] - \phi(H_{\{f,g\}},\,H_h,\,\cdot\,)\notag \\ 
    &\quad + \bigl[[dg,\,dh],\,df\bigr] + \bigl[\phi(H_g,\,H_h,\,\cdot\,),\,df\bigr] 
                                              - \phi(H_{\{g,h\}},\,H_f,\,\cdot\,)\notag \\ 
    &\quad + \bigl[[dh,\,df],\,dg\bigr] + \bigl[\phi(H_h,\,H_f,\,\cdot\,),\,dg\bigr] 
                - \phi(H_{\{h,f\}},\,H_g,\,\cdot\,) + d\bigl(\phi(H_f,\,H_g,\,\cdot\,)\bigr).
\end{align}
Note that $[df,\,dg] = -\{df,\,dg\} - \phi(H_f,\,H_g,\,\cdot\,)$ is used in the above calculation. 
For any $f,\,g\in C^\infty(P)$, we define $\eta_{fg}\in \varGamma(T^*P)$ as $\eta_{fg}(X):= \phi(H_f,\,H_g,\,X)$. 
Then, we have 
\begin{align}\label{algebroid2}
 [\eta_{fg},\,dh](H_k)\nonumber
 =&\, -i_{H_h}d\eta_{fg}(H_k) - d\bigl(\Pi(\eta_{fg},\,dh)\bigr)(H_k) - \phi(\sharp\eta_{fg},\,H_h,\,H_k)\nonumber\\
 =&\, \eta_{fg}(H_k,\,H_h) - d\bigl(\eta_{fg}(H_h)\bigr)(H_k) - \phi(\sharp\eta_{fg},\,H_h,\,H_k)\nonumber\\ 
 =&\, H_k\bigl(\eta_{fg}(H_h)\bigr) - H_h\bigl(\eta_{fg}(H_k)\bigr) + \eta_{fg}([H_f,\,H_g])\nonumber\\
        &\qquad\qquad\qquad\qquad -d\bigl(\phi(H_f,\,H_g,\,H_h)\bigr)(H_k) - \phi(\sharp\eta_{fg},\,H_h,\,H_k)\notag\\
 =&\, H_k\bigl(\phi(H_f,\,H_g,\,H_h)\bigr) - H_h\bigl(\phi(H_f,\,H_g,\,H_k)\bigr) 
                                              + \phi\bigl(H_f,\,H_g,\,[H_h,\,H_k]\bigr)\notag\\ 
      &\qquad\qquad\qquad\qquad - d\bigl(\phi(H_f,\,H_g,\,H_h)\bigr)(H_k) - \phi(\sharp\eta_{fg},\,H_h,\,H_k),
\end{align}
for any $k\in C^\infty(P)$. From Proposition \ref{eq2}~, 
\[
 dh\bigl(H_{\{f,g\}}+[H_f,\,H_g]\bigr) = \phi(H_f,\,H_g,\,H_h) = \eta_{fg}\bigl(\sharp(dh)\bigr) = -dh(\sharp\eta_{fg}).
\]
Thus, 
\begin{equation}\label{algebroid3}
 \phi(H_{\{f,g\}},\,H_h,\,H_k) = -\phi([H_f,\,H_g],\,H_h,\,H_k) - \phi(\sharp\eta_{fg},\,H_h,\,H_k).
\end{equation}
It follows from (\ref{algebroid2})~and (\ref{algebroid3})~that 
\begin{align*}
 (\ref{algebroid1}) 
=&\, \bigl[[df,\,dg],\,dh\bigr] + \bigl[[dg,\,dh],\,dk\bigr] + \bigl[[dh,\,df],\,dg\bigr]\\
 &\: + H_k\bigl(\phi(H_f,\,H_g,\,H_h)\bigr) - H_h\bigl(\phi(H_f,\,H_g,\,H_k)\bigr) 
                                               + \phi\bigl(H_f,\,H_g,\,[H_h,\,H_k]\bigr)\\ 
 &\;- d\bigl(\phi(H_f,\,H_g,\,H_h)\bigr)(H_k) - \phi(\sharp\eta_{fg},\,H_h,\,H_k) 
                                                      + H_k\bigl(\phi(H_g,\,H_h,\,H_f)\bigr)\\
 &\; - H_f\bigl(\phi(H_g,\,H_h,\,H_k)\bigr) + \phi\bigl(H_g,\,H_h,\,[H_f,\,H_k]\bigr) 
                                                   - d\bigl(\phi(H_g,\,H_h,\,H_f)\bigr)(H_k)\\ 
 &\; - \phi(\sharp\eta_{gh},\,H_f,\,H_k) + H_k\bigl(\phi(H_h,\,H_f,\,H_g)\bigr) - H_g\bigl(\phi(H_h,\,H_f,\,H_k)\bigr)\\ 
 &\; + \phi\bigl(H_h,\,H_f,\,[H_g,\,H_k]\bigr) - d\bigl(\phi(H_h,\,H_f,\,H_g)\bigr)(H_k) - \phi(\sharp\eta_{hf},\,H_g,\,H_k)\\ 
 &\; + \phi([H_f,\,H_g],\,H_h,\,H_k) + \phi(\sharp\eta_{fg},\,H_h,\,H_k) + \phi([H_g,\,H_h],\,H_f,\,H_k)\\ 
 &\; + \phi(\sharp\eta_{gh},\,H_f,\,H_k) + \phi([H_h,\,H_f],\,H_g,\,H_k) + \phi(\sharp\eta_{hf},\,H_g,\,H_k)\\ 
=&\; \bigl[[df,\,dg],\,dh\bigr] + \bigl[[dg,\,dh],\,dk\bigr] + \bigl[[dh,\,df],\,dg\bigr]\\
 &\; - H_f\bigl(\phi(H_g,\,H_h,\,H_k)\bigr) - H_g\bigl(\phi(H_h,\,H_f,\,H_k) - H_h\bigl(\phi(H_f,\,H_g,\,H_k)\bigr)\\
 &\; + H_k\bigl(\phi(H_f,\,H_g,\,H_h)\bigr) + \phi([H_f,\,H_g],\,H_h,\,H_k) + \phi([H_h,\,H_f],\,H_g,\,H_k)\\
 &\; + \phi\bigl(H_g,\,H_h,\,[H_f,\,H_k]\bigr) + \phi\bigl(H_g,\,H_h,\,[H_f,\,H_k]\bigr) 
                                                            + \phi\bigl(H_h,\,H_f,\,[H_g,\,H_k]\bigr)\\
 &\; + \phi\bigl(H_f,\,H_g,\,[H_h,\,H_k]\bigr) + 2H_k\bigl(\phi(H_g,\,H_h,\,H_f)\bigr) 
                                                          - 2d\bigl(\phi(H_g,\,H_h,\,H_f)\bigr)(H_k)\\
=& \bigl[[df,\,dg],\,dh\bigr] + \bigl[[dg,\,dh],\,dk\bigr] + \bigl[[dh,\,df],\,dg\bigr] - (d\phi)(H_f,\,H_g,\,H_h,\,H_k)\\
=&\, \bigl[[df,\,dg],\,dh\bigr] + \bigl[[dg,\,dh],\,dk\bigr] + \bigl[[dh,\,df],\,dg\bigr]. 
\end{align*}
Consequently, $[\cdot,\,\cdot]_\phi$ satisfies the Jacobi identity. 

\begin{dfn}
 Let $A\to M$ be a Lie algebroid with anchor map $\sharp:A\to TM$. 
  A left(right) action of $A$ on a smooth manifold $N$ consists of a smooth map 
 $J:N\to M$ called the moment map, 
and a Lie algebra (anti)homomorphism $\varrho_N:\varGamma(A)\to \varGamma(TN)$ which satisfy the following 
 conditions{\rm :}
 \begin{enumerate}[{\quad \rm(1)}]
   \item $dJ\circ \varrho_N(\alpha) \,=\, \sharp\alpha;$
   \item $\varrho_N(f\alpha) \,= \, (J^*f)\varrho_N(\alpha),$
 \end{enumerate}
 for any $f\in C^\infty(M)$ and $\alpha\in \varGamma(A)$.
\end{dfn}

If $(P,\,\Pi_P,\,\phi_P),\,(Q,\,\Pi_Q,\,\phi_Q)$ are twisted Poisson manifolds, then any twisted Poisson map 
$J:Q\to P$ induces a Lie algebroid action of $T^*P$ on $Q$ by 
\[\varrho_Q: \varGamma(T^*P)\longrightarrow \varGamma(TQ),\quad \alpha \longmapsto \Pi_Q^\sharp(J^*\alpha).\]
In fact, by (\ref{twPoisson map})~ and Proposition \ref{twPmap}, 
        for any $f\in C^{\infty}(Q)$ and $\alpha,\,\beta\in \varGamma(T^*P)$, 
\begin{align*}
 &\Bigl(\Pi_Q^\sharp\bigl([J^*\alpha,\,J^*\beta]_{\phi_Q}\bigr)\Bigr)f\\ 
 =\,& -[J^*\alpha,\,J^*\beta]_{\phi_Q}(H_f)\\
 =\,& -\Bigl(i_{\Pi_Q^\sharp(J^*\alpha)}d(J^*\beta) - i_{\Pi_Q^\sharp(J^*\beta)}d(J^*\alpha) 
              - d\bigl(\Pi_Q(J^*\alpha,\,J^*\beta)\bigr) \,-\, 
          \phi_Q\bigl(\Pi_Q^\sharp(J^*\alpha),\,\Pi_Q^\sharp(J^*\beta),\,\cdot\bigr)\Bigr)(H_f)\\
 =\,& -J^*(d\beta)\bigl(\Pi_Q^\sharp(J^*\alpha),\,H_f\bigr) + J^*(d\alpha)\bigl(\Pi_Q^\sharp(J^*\beta),\,H_f\bigr)\\
     &\qquad\qquad\qquad\qquad\qquad\qquad\qquad+ J^*\bigl(d\Pi_Q(\alpha,\,\beta)\bigr)(H_f) \,+\, 
          \phi_Q\bigl(\Pi_Q^\sharp(J^*\alpha),\,\Pi_Q^\sharp(J^*\beta),\,H_f\bigr)\\
 =\,& -d\beta(\Pi_P(\alpha),\,(dJ)(H_f)) + d\alpha(\Pi_P(\beta),\,(dJ)(H_f))\\
       &\qquad\qquad\qquad\qquad\qquad\qquad\qquad+ d\bigl(\Pi_P(\alpha,\,\beta)\bigr)((dJ)(H_f)) \,+\, 
             (J^*\phi_P)\bigl(\Pi_Q^\sharp(J^*\alpha),\,\Pi_Q^\sharp(J^*\beta),\,H_f\bigr)\\
 =\,& -d\beta(\Pi_P(\alpha),\,(dJ)(H_f)) + d\alpha(\Pi_P(\beta),\,(dJ)(H_f))\\ 
  &\qquad\qquad\qquad\qquad\qquad\qquad\qquad + d\bigl(\Pi_P(\alpha,\,\beta)\bigr)((dJ)(H_f))
                                  + \phi_P\bigl(\Pi_P^\sharp(\alpha),\,\Pi_P^\sharp(\beta),\,(dJ)(H_f)\bigr)\\
 =\, & -\Bigl(i_{\Pi_P^\sharp(\alpha)}d\beta - i_{\Pi_P^\sharp}(\beta)d\alpha
          - d\bigl(\Pi_P(\alpha,\,\beta)\bigr) - 
             \phi_P\bigl(\Pi_P^\sharp(\alpha),\,\Pi_P^\sharp(\beta),\,\cdot\bigr)\Bigr)(H_f)\\
 =\, & -\bigl(J^*[\alpha,\,\beta]_{\phi_P}\bigr)(H_f)\\
 =\, & \bigl(\Pi_Q^\sharp\bigl(J^*[\alpha,\,\beta]_{\phi_P}\bigr)\bigr)f
\end{align*}
 
It follows from this that $\varrho_Q\bigl([\alpha,\,\beta]_{\phi_Q}\bigr) = [\varrho_Q(\alpha),\,\varrho_Q(\beta)]_{\phi_P}$.  
In other words, $\varrho_Q$ is a Lie algebra homomorphism. 
Moreover, it can be verified that, for any $f\in C^\infty(P)$ and $\alpha\in \varGamma(T^*P)$, 
\begin{align*}
 \varrho_Q(f\alpha) &= \Pi_Q^\sharp(J^*(f\alpha))=J^*f\Pi_Q(J^*\alpha)=(J^*f)\varrho_Q(\alpha),\\
 dJ\bigl(\varrho_Q(\alpha)\bigr) &= dJ\Pi_Q^\sharp(J^*\alpha) = \Pi_P^\sharp(\alpha).
\end{align*}
This leads us to the conclusion that any twisted Poisson map $J:Q\to P$ is equipped with a Lie algebroid action of $T^*P$.

\vspace*{0.6cm}
As is well known, given a Lie groupoid $\varGamma\rightrightarrows M$, one can construct the Lie algebroid over $M$ 
denoted by $\mathcal{A}(\varGamma)$. For a full discussion of the Lie algebroid of the Lie groupoid, we refer to 
Crainic, M. and Fernandes, R.-L. \cite{cf2}. A Lie algebroid $A\to M$ is said to be integrable 
if there exists a Lie groupoid $\varGamma\rightrightarrows M$ 
such that $\mathcal{A}(\varGamma)$ is isomorphic to $A$. 

\begin{dfn}
A twisted Poisson manifold is said to be integrable if its cotangent bundle is integrable as Lie algebroid. 
\end{dfn}

The integrability problem of Lie algebroids was studied by many people, for instance, Pradines, J.~\cite{pradines},\,
Mackenzie, K.~\cite{mackenzie} and Crainic, M. and Fernandes, R.-L.~\cite{cf}.  
The solution of integrability problem of twisted Poisson manifolds was given by Cattaneo, A. and Xu, P.~(\cite{cx}). 
They proved the following result:

\begin{thm}
{\rm (Cattaneo, A. and Xu, P.)}~There is a bijection between integrable twisted Poisson structures and twisted symplectic 
groupoids which are source-simply connected. 
\end{thm}

That is, twisted Poisson manifolds may be integrated to twisted symplectic groupoids. For an integrable twisted Poisson 
manifold $P$, we denote by $\mathcal{G}(P)$ the twisted symplectic groupoid associated with $P$ in the above theorem. 
We refer to Definition \ref{twsymplectic} for a twisted symplectic groupoid. 

\section{Geometric Morita equivalence}
\subsection{Morita invariants}
First, we will introduce the notion of Morita equivalence of twisted Poisson manifolds and 
exhibit some examples. 

\begin{dfn}\label{dfn:morita}{\rm (\cite{yuji1},\cite{yuji2})}
Let $P_i$ be integrable $\phi_i$-twisted Poisson manifolds~(i=1,2). $P_1$ and $P_2$ are 
said to be {\rm (}strong{\rm )} Morita equivalent if there exist a smooth manifold $S$ equipped with a 
non-degenerate 2-form $\omega_S$ and surjective submersions $J_i:S\to P_i$ such that

 \begin{enumerate}[{\quad \rm(1)}]
  \item $(S,\,\omega_S)$ is a $(J_1^*\phi_1-J_2^*\phi_2)$-twisted symplectic manifold{\rm ;}
  \item $J_1$ is a complete t.s.realization, and $J_2$ is a complete anti-t.s.realization{\rm ;}
  \item Each $J_i$-fiber~$(i=1,2)$ is connected, and simply-connected{\rm ;}
  \item The subspaces $\ker(dJ_1)_x,\,\ker(dJ_2)_x$ of $T_xS~ (\forall x\in S)$ are symplectically 
   orthogonal to one another{\rm :}
   \[\bigl(\ker(dJ_1)_x\bigr)^\perp = \ker(dJ_2)_x\quad \text{and} \quad \bigl(\ker(dJ_2)_x\bigr)^\perp = \ker(dJ_1)_x,\]
 where
   \[\bigl(\ker(dJ_i)_x\bigr)^\perp = \bigl\{\,\boldsymbol{u}\in T_xS \bigm| 
       \omega_S(\boldsymbol{u},\,\boldsymbol{v})=\boldsymbol{0}~
                 \bigl(\forall \boldsymbol{v}\in \ker(dJ_i)_x\bigr)\,\bigr\},\quad (i=1,2).\]
 \end{enumerate}
\end{dfn}

A twisted symplectic manifold $S$ in Definition \ref{dfn:morita} is called a 
$(P_1,\,P_2)$-equivalence bimodule (or an equivalence bimodule for short), 
and denoted by $P_1\overset{J_1}{\leftarrow}S\overset{J_2}{\rightarrow}P_2$. 

\begin{ex}{\rm 
An integrable twisted Poisson manifold is Morita equivalent to itself with an equivalence bimodule $\mathcal{G}(P)$.
}
\end{ex}

\begin{ex}{\rm
$\rm (Example 2.1 in \cite{xu2})$~ Let $S$ be a connected and simply-connected symplectic manifold, and $M$ a connected 
smooth manifold with a trivial Poisson structure${\rm :}\{\cdot,\,\cdot\}\equiv 0$. 
Then, $S\times M$ is Morita equivalent to $M$ with a equivalence bimodule $S\times T^*M$. 
}
\end{ex}

\begin{ex}{\rm 
Let $P_i$ and $Q_i$ be twisted Poisson manifolds~$(i=1,2)$. 
Assume that $P_1$ and $Q_1$ are Morita equivalent to $P_2$ and $Q_2$ respectively, with equivalence bimodules 
$P_1\overset{J_1}{\leftarrow}X\overset{J_2}{\rightarrow}P_2$ 
and $Q_1\overset{J'_1}{\leftarrow}Y\overset{J'_2}{\rightarrow}Q_2$. 
Then, $P_1\times Q_1$ and $P_2\times Q_2$ are Morita equivalent${\rm :}~
P_1\times Q_1\overset{J_1\times J_1'}{\longleftarrow}X\times Y\overset{J_2\times J_2'}{\longrightarrow}P_2\times Q_2$. 
}
\end{ex}

\begin{ex}{\rm 
Two simply-connected twisted symplectic manifolds $(S_i,\,\omega_i,\,\psi_i)~(i=1,2)$ are Morita equivalent 
each other. In fact, 
we denote the natural projections from $S_1\times S_2$ to $S_i$ by ${\rm pr}_i~(i=1,2)$ and 
set $\omega = {\rm pr}_1^*\omega_1 - {\rm pr}_2^*\omega_2$. It is easy to verify that ${\rm pr}_1:S_1\times S_2\to S_1$ 
and ${\rm pr}_2:S_1\times S_2\to S_2$ are t.s.realization and anti-t.s.realization, respectively. 
If $\boldsymbol{V}=V_1\oplus V_2\in \bigl(\ker(d{\rm pr}_1)_x\bigr)^\perp~(x\in S_1\times S_2)$ satisfies that 
$\omega(\boldsymbol{V},\,\boldsymbol{U})=0$ for any $\boldsymbol{U}\in \ker(d{\rm pr}_1)_x$, 
then the non-degeneracy of $\omega_2$ implies that $\boldsymbol{V}\in \ker(d{\rm pr}_2)_x$. 
From this, we have $\bigl(\ker(d{\rm pr}_1)_x\bigr)^\perp = \ker(d{\rm pr}_2)_x$. Similarly, it is also shown that 
$\bigl(\ker(d{\rm pr}_2)_x\bigr)^\perp = \ker(d{\rm pr}_1)_x$. Therefore, $S_1$ and $S_2$ are Morita equivalent with a 
equivalence bimodule $S_1\times S_2$. 
}
\end{ex}

Morita equivalence is indeed an equivalence relation among twisted Poisson manifolds: As for the transitivity, 
suppose that $(S_1,\,\omega_1)$ is a $(P_1,\,P_2)$-equivalence bimodule with moment maps 
$P_1\overset{J_1}{\leftarrow}S_1\overset{J_2}{\to}P_2$ 
and $(S_2,\,\omega_2)$ is a $(P_2,\,P_3)$-equivalence bimodule with moment maps 
$P_2\overset{J_2'}{\leftarrow}S_2\overset{J_3}{\to}P_3$. We define a smooth manifold $S_1\otimes S_2$ to be the 
quotient of the fiber product by its characteristic foliation:
\begin{equation}\label{tensor}
S_1\otimes S_2 = (S_1\times_{P_2} ^{J_2,J'_2}S_2)\big/\ker \bigl(\iota^*(\omega_1\oplus\omega_2)\bigr),
\end{equation} 
where $\iota:S_1\times_{P_2} ^{J_2,J'_2}S_2\hookrightarrow S_1\times S_2$ is the canonical embedding map. 
In addition, we define a non-degenerate 2-form $\tilde{\omega}$ on $S_1\otimes S_2$ by 
\[\tilde{\omega}_{\pi(p)}(\pi_*u,\,\pi_*v) = (\omega_1\oplus\omega_2)_p(\iota_*u,\,\iota_*v),\]
where $\pi:S_1\times_{P_2} ^{J_2,J'_2}S_2\to S_1\otimes S_2$ is the natural projection. 
It is verified that the 2-form $\tilde{\omega}$ is well-defined in a way similar to \cite{xu3}.
Then $(S_1\otimes S_2,\,\tilde{\omega})$ is a $(P_1,\,P_3)$-equivalence bimodule with the moment maps
$\tilde{J}_1:S_1\otimes S_2\to P_1$ and $\tilde{J}_3:S_1\otimes S_2\to P_3$ given by 
$\tilde{J}_1([x,\,y]):=J_1(x)$ and $\tilde{J}_3([x,\,y]):=J_3(y)$, respectively.   

\vspace*{0.5cm}
\noindent {\bf Remark}.~
We note that the tensor product in (\ref{tensor})~is not associative, but just associative up to a 
bimodule isomorphism. For any integrable twisted Poisson manifolds $P_1$ and $P_2$, we denote 
by $\mathcal{B}(P_1,\,P_2)$ the set of all $(P_1,\,P_2)$-equivalence bimodules. It is verified that 
$\mathcal{B}(P_1,\,P_2)$ forms a category whose morphisms are complete twisted Poisson maps $f$ between 
equivalence bimodules $P_1\overset{J_1}{\leftarrow}S_1\overset{J_2}{\to}P_2$ and 
$P_1\overset{K_1}{\leftarrow}S_2\overset{K_2}{\to}P_2$ which satisfies $J_1=K_1\circ f$ and $J_2=K_2\circ f$. 
Then, we have the bicategory $\mathsf{twPoiss}$ which has integrable twisted Poisson manifolds as 0-cells, 
equivalence bimodules as 1-cells and the tensor product as compositions. Two integrable twisted Poisson manifolds 
are Morita equivalent if and only if they are isomorphic object in $\mathsf{twPoiss}$. For the definition of a bicategory, 
we refer to \cite{benabou}. 

\vspace*{0.5cm}
As in the case of Poisson manifolds, twisted Poisson manifolds which are Morita equivalent have similar features. 
For instance, they have isomorphic first cohomology groups and representation categories. For a twisted symplectic 
manifold, its fundamental group is a complete Morita invariant. 
We revise these results discussed in \cite{yuji1}. 
\begin{thm}
Let $(S_i,\,\omega_i,\,\psi_i)~(i=1,2)$ be connected and integrable twisted symplectic manifolds. 
$S_1$ and $S_2$ are Morita equivalent if and only if their fundamental groups are isomorphic each other. 
\end{thm}

\noindent ({\it Proof})~
Assume that $S_1$ and $S_2$ are Morita equivalent with a $(S_1,\,S_2)$-equivalence bimodule 
$S_1\overset{J_1}{\leftarrow} S\overset{J_2}{\rightarrow}S_2$. By the assumption that $J_1:S\to S_1$ is 
complete and $J_1$-fiber is connected and simply-connected, we have the following exact sequence of homotopy groups:
\[\{0\}\cong \pi_1\bigl(J_1^{-1}(*)\bigr)\rightarrow \pi_1(S)\rightarrow \pi_1(S_1)\rightarrow 
\pi_0\bigl(J_1^{-1}(*)\bigr)\cong \{0\},\]
where $\ast$ denotes any fixed point in $S$. This implies that $\pi_1(S_1)\cong \pi_1(S)$. 
It is also verified that $\pi_1(S_2)\cong \pi_1(S)$ in a similar way. 
Therefore, $\pi_1(S_1)$ is isomorphic to $\pi_1(S_2)$. 

Conversely, suppose that the fundamental groups of $S_1$ and $S_2$ are isomorphic. 
We take the universal covering spaces of $S_i~(i=1,2)$ and denote them by $\widetilde{S}_i$. 
Then, $G\cong \pi_1(S_i)$ acts on $\widetilde{S}_1\times \widetilde{S}_2$ diagonally. 
We define the maps $\widetilde{\rho}_i~(i=1,2)$ from the quotient space $S=(\widetilde{S}_1\times\widetilde{S}_2)/G$ to 
$S_i\cong\widetilde{S}_i/G$ by 
\[\widetilde{\rho}_i\bigl([x_1,\,x_2]\bigr) \,=\, \rho_i(x_i),\quad \bigl(\forall [x_1,\,x_2]\in S\bigr),\]
where $\rho_i:\widetilde{S}_i\to S_i$ denotes the covering map of $S_i$.  
Using these maps, we define a $2$-form $\omega_S$ on $S$ 
by $\omega_S := \widetilde{\rho}_1^*\omega_1 - \tilde{\rho}_2^*\omega_2$. 
It is easy to checked that $\omega_S$ is non-degenerate and 
$d\omega_S = \widetilde{\rho}_1^*\psi_1 - \widetilde{\rho}_2^*\psi_2$, 
that is, $(S,\,\omega_S)$ is a $(\widetilde{\rho}_1^*\psi_1 - \widetilde{\rho}_2^*\psi_2)$-twisted symplectic manifold. 
Let $\{\cdot,\,\cdot\}_S,\,\{\cdot,\,\cdot\}_1$ and $\{\cdot,\,\cdot\}$ be the brackets induced from 
$\omega_S,\,\omega_1$ and $\rho_1^*\omega_1-\rho_2^*\omega_2$, respectively. Then, for any $f,\,g\in C^{\infty}(S_1)$, 
\begin{align*}
\{\widetilde{\rho}^*_1f,\,\widetilde{\rho}^*_2g\}_S\bigl([x_1,\,x_2]\bigr) 
        &= \bigl\{\pi^*(\widetilde{\rho}_1^*f),\,\pi^*(\widetilde{\rho}_2^*g)\bigr\}(x_1,\,x_2)\\
    &= \{(\widetilde{\rho}_1\circ\pi)^*f(\,\cdot,\, x_2),\,(\widetilde{\rho}_1\circ\pi)^*g(\,\cdot,\, x_2)\}_1(x_1)\\
    &= \{\rho_1^*f,\,\rho_1^*g\}(x_1) = \{f,\,g\}_1\bigl(\rho_1(x_1)\bigr)\\
    &= \bigl(\widetilde{\rho}_1^*\{f,\,g\}_1\bigr)\bigl([x_1,\,x_2]\bigr).
\end{align*}

In addition, the Hamiltonian vector field $H_{\widetilde{\rho}_1^*f}=H_{\rho^*_1f}$ is also complete if $H_f$ is complete. 
These imply that $\widetilde{\rho}_1:S\to S_1$ is a complete t.s.realization. Similarly, 
$\widetilde{\rho}_2:S\to S_2$ is a complete anti-t.s.realization. 
It is also verified that $\bigl(\ker (d\widetilde{\rho}_1)_x\bigr)^\perp = \ker (d\widetilde{\rho}_2)_x$ and 
$\bigl(\ker (d\widetilde{\rho}_2)_x\bigr)^\perp = \ker (d\widetilde{\rho}_1)_x$ for any $x\in S$.
Accordingly, $S_1$ and $S_2$ are Morita equivalent. \qquad\qquad\qquad\qquad\qquad\qquad\qquad\quad $\Box$

\vspace{0.5cm}
Let $(P,\,\Pi,\,\phi)$ be a twisted Poisson manifold. 
As is discussed in 1.2, the space of 1-forms $\varGamma(T^*P)$ carries Lie algebra structure. 
Using the Lie bracket $[\cdot,\,\cdot]_\phi$, we define the differential operator on the space of multi-vector 
fields on $P$ by 
\begin{align*}
 (d_{\Pi,\phi} A)(\alpha_1,\cdots,\alpha_{k+1})\,:=&\, -\sum _{i=1}^{k+1}(-1)^{i+1}\Pi^\sharp(\alpha_i)
                                         \bigl( A(\alpha_1,\cdots,\widehat{\alpha_i},\cdots,\alpha_{k+1})\bigr)\\
 &\,-\, \sum_{i<j}(-1)^{i+j}A\bigl( [\alpha_i,\, \alpha_j]_\phi,\,\alpha_1,\cdots,\widehat{\alpha_i},
                  \cdots,\widehat{\alpha_j},\cdots,\alpha_{k+1}\bigr), 
\end{align*}
for any $\alpha_1,\cdots,\alpha_{k+1}\in \varGamma(T^*P)$. Twisted Poisson cohomology is the cohomology of 
the differential complex 
$\bigl(\varGamma(\wedge^\bullet T^*P),\,d_{\Pi,\phi}\bigr)$, denoted by $H^\bullet_{\Pi,\phi}(P)$
(see Kosmann-Schwarzbach, Y. and Laurent-Gengoux, C. \cite{kosmann} and \cite{severa}). 
We remark that the following formulae hold:
\begin{align*}
 &d_{\Pi,\,\phi}f = [\Pi,\,f] = H_f,\\ 
 &(d_{\Pi,\,\phi}X)(df,\,dg) = -(\mathcal{L}_X\Pi)(df,\,dg) - \Phi(H_f,\,H_g,\,X).
\end{align*}

\begin{thm}
Let $(P_1,\,\Pi_1,\,\phi_1)$ and $(P_2,\,\Pi_2,\,\phi_2)$ be integrable twisted Poisson manifolds. 
If they are Morita equivalent, then their first cohomology groups $H_{\Pi_1}^1(P_1)$ and $H_{\Pi_2}^1(P_2)$ are isomorphic. 
\end{thm}
\noindent ({\it Proof})~
Assume that $P_1$ and $P_2$ are Morita equivalent with a equivalent bimodule 
$P_1\overset{J_1}{\leftarrow}S\overset{J_2}{\rightarrow}P_2$. Let $X$ be a vector field on $P_1$ such that 
$d_{\Pi_1,\,\phi_1}X=0$. Since $J_1:S\to P_1$ is a surjective submersion and 
$\{J_1^*C^{\infty}(P_1),\,J_2^*C^{\infty}(P_2)\}_S=0$, any tangent vector to the $J_2$-fiber at $x\in S$ is represented by 
the value of Hamiltonian vector field of a function in $J_1^*C^\infty(P_1)$. On the basis of this observation, we define 
a 1-form $\theta_X$ on each $J_2$-fiber by 
\[\langle (\theta _X)_x,\, \varphi(x)(H _{J_1^{*}f})_x\rangle = \varphi(x)(Xf)_{J_1(x)},\quad \bigl(x\in S,\, 
\varphi\in C^{\infty}(S),\,f\in C^{\infty}(P_1)\bigr),\]
where $H_\bullet$ means the Hamiltonian vector field with respect to the $(J_1^*\phi_1 - J_2^*\phi_2)$-twisted 
symplectic form. Note that $X_{J_1(x)}$ is represented by 
$X_{J_1(x)}=(dJ_1)_x(H_{J_1^*h})$ for some $h\in C^\infty(P_1)$. 
If $(H_{J_1^*f})_x=(H_{J_1^*g})_x$, then we have  
\begin{align*}
\langle (\theta _X)_x,\, (H _{J_1^{*}f})_x\rangle 
  &= (Xf)_{J_1(x)} = (df)_{J_1(x)}(X_{J_1(x)})=d(J_1^*f)_x\bigl((H_{J_1^*h})_x\bigr)\\
  &= -d(J_1^*h)_x\bigl((H_{J_1^*f})_x\bigr) = -d(J_1^*h)_x\bigl((H_{J_1^*g})_x\bigr)\\
  &= \langle (\theta _X)_x,\, (H _{J_1^{*}g})_x\rangle .
\end{align*}
Hence, the 1-form $\theta_X$ is well-defined. 
As will be shown in Lemma \ref{claim}, $\theta_X$ is closed on each $J_2$-fiber. 
Assuming this claim for the moment, we continue the proof. 

Since each $J_2$-fiber is closed and simply-connected, $\theta_X$ is exact thereon. As mentioned in \cite{lu}, 
there exists a smooth function $\varLambda\in C^\infty(S)$ such that $\theta_X=d\varLambda$ and 
\begin{equation}\label{cohomology1}
 H_{\varLambda}\!\cdot\!\bigl(J_1^{*}C^{\infty}(P_1)\bigr)\subset J_1^{*}C^{\infty}(P_1)
\end{equation}
when restricted to each $J_2$-fiber. 

Now, for any $f\in C^\infty(P_1)$ and $g\in C^\infty(P_2)$
\begin{align*}
 &\left\{H_{\varLambda}J_2^{*}g,\, J_1^{*}f\right\}_S\\ =& \left\{\{J_2^{*}g,\, \varLambda \}_S, J_1^{*}f \right\}_S\\
  =&\, -\bigl\{\{\varLambda , J_1^{*}f\}_S,\, J_2^{*}g\bigr\}_S - \bigl\{\{J_1^{*}f,\, J_2^{*}g\}_S,\varLambda \bigr\}_S
                    - (J_1^{*}\phi _1-J_2^{*}\phi _2)(H_{J_2^{*}g},H_{\varLambda},H_{J_1^{*}f})\\
  =&\, \bigl\{H_\varLambda(J_1^{*}f),\, J_2^{*}g\bigr\}_S - \bigl\{\{J_1^{*}f,\, J_2^{*}g\}_S,\varLambda \bigr\}_S
                              - (J_1^{*}\phi _1-J_2^{*}\phi _2)(H_{J_2^{**}g}, H_{\varLambda}, H_{J_1^{*}f})\\ 
  =&\,\bigl\{H_\varLambda(J_1^{*}f),\, J_2^{*}g\bigr\}_S - \bigl\{\{J_1^{*}f,\, J_2^{*}g\}_S,\varLambda \bigr\}_S\\
             &\qquad\qquad- \phi _1\bigl({J_1}_*H_{J_2^{*}g},\, {J_1}_*H_{\varLambda},\, {J_1}_*H_{J_1^{*}f}\bigr) + 
                     \phi _2\bigl({J_2}_*H_{J_2^{*}g},\, {J_2}_*H_{\varLambda},\, {J_2}_*H_{J_1^{*}f}\bigr) 
\end{align*}
From (\ref{cohomology1})~and assumption, the last two terms of the above formula are equal to zero.
As a result, we have that 
\[H_{\varLambda}\bigl(J_2^*C^{\infty}(P_2)\bigr)\subset J_2^*C^{\infty}(P_2).\]
Therefore, the push-forward vector field on $P_2$ given by $X'= (J_2)_*H_\varLambda$ is well-defined. For any 
$f,\,g\in C^\infty(P_2)$, 
\begin{align*} 
 &\phi_2(\Pi_2^\sharp(df),\,\Pi_2^\sharp(dg),\, X')\\ =&\, -\bigl(d\Pi_2(df,\,dg)\bigr)(X') - 
                                      \Pi_2\bigl(d(X'g),\,df\bigr) + \Pi_2\bigl(d(X'f),\,dg\bigr)\\
 =&\, -X'(\Pi_2(df,\,dg)) - \Pi_2(d(X'g),\,df) + \Pi_2(d(X'f),\,dg)\\
 =&\, -X'\left( \Pi _2(df,\, dg)\right) + \Pi _2(df,\, \mathcal{L}_{X'}dg) + \Pi _2(\mathcal{L}_{X'}df,\, dg)\\
 =&\, -(\mathcal{L}_{X'}\pi _2)(df,\,dg).
\end{align*}
This implies that $d_{\Pi_2,\phi_2}X' = 0$. The cohomology class of $X'$ does not depend on the choice of $\varLambda$ 
and the map which corresponds $X$ to the cohomology class of $X'$ induces a map from $H_{\phi_1}^1(P_1)$ 
to $H_{\phi_2}^1(P_2)$~(see \cite{lu}). 
We can construct a map from $H_{\phi_2}^1(P_2)$ to $H_{\phi_1}^1(P_1)$ in a similar way. 
Hence, $H_{\phi_1}^1(P_1)$ and $H_{\phi_2}^1(P_2)$ are isomorphic each other. 
\qquad\qquad\qquad\qquad\qquad\qquad\qquad\qquad\qquad$\quad \Box$

\vspace*{0.5cm}
It remains us to show the following claim:
\begin{lem}\label{claim}
The 1-form $\theta_X$ is closed on each $J_2$-fiber.
\end{lem}
\noindent ({\it Proof})~
Note that $J_1\times J_2:S\to P_1\times \overline{P_2},\,p\mapsto \bigl(J_1(p),\,J_2(p)\bigr)$ is a twisted Poisson 
map. For any $f,\,g\in C^\infty(P_1)$, we have 
\begin{align*} 
(d\theta _X)(H_{J_1^*f},\, \xi _{J_1^*g}) &= H_{J_1^*f}\left( \theta _X(H_{J_1^*g})\right ) - 
                      H_{J_1^*g}\left( \theta _X(H_{J_1^*f})\right ) - \theta _X\left([H_{J_1^*f},\, H_{J_1^*g}]\right)\\
 &= H_{J_1^*f}\left( J_1^*(Xg)\right) - 
                      H_{J_1^*g}\left( J_1^*(Xf)\right ) - \theta _X\left([H_{J_1^*f},\, H_{J_1^*g}]\right)\\
 &= \{J_1^*(Xg),\,J_1^*f\}_S - \{J_1^*(Xf),\, J_1^*g\}_S - \theta_X\left([H _{J_1^*f},\, H _{J_1^*g}]\right)\\ 
 &= J_1^*\{Xg,\, f\} _1 + J_1^*\{g,\, Xf\} _1 - \theta _X\left([H_{J_1^*f},\, H_{J_1^*g}]\right), 
\end{align*}
where $\{\cdot,\, \cdot\}_1$ means the bracket induced from the bivector $\Pi_1$. 
Given $f,\,g\in C^\infty(P_1)$, we define the map between the smooth functions on each $J_2$-fiber by 
\begin{equation}
F\in C^{\infty}(S) \longmapsto J_1^*\phi _1(H_{J_1^*f},\, H_{J_1^*g},\, H_F) \in C^{\infty}(S).
\end{equation}
This map is considered as a vector field on each $J_2$-fiber, so, can be represented by $H_{J_1^*k}$ for 
some $k\in C^\infty(P_1)$. Accordingly, using by Proposition \ref{eq2}, we have 
\[\theta_X\bigl( [H_{J_1^*f},\, H_{J_1^*g}] + H_{\{J_1^*f,\,J_1^*g\}}\bigr) = \theta_X(H_{J_1^*k}).\]
Therefore, from assumption it follows that
\begin{align*}\label{closed}
(d\theta _X)(H_{J_1^*f},\, H_{J_1^*g}) &= J_1^*\{Xg,\, f\} _1 + J_1^*\{g,\, Xf\} _1 
                                      + \theta_X\bigl( H_{J_1^*\{f,\,g\}_1}\bigr) - \theta _X(H_{J_1^*k})\\
                  &= J_1^*\{Xg,\, f\} _1 + J_1^*\{g,\, Xf\} _1 + J_1^*X\{f,g\}_1 - \theta _X(H_{J_1^*k})\\                                    
                  &= J_1^*\left( \{Xg,\, f\} _1 + \{g,\, Xf\} _1 + X\{f,g\}_1\right) - \theta _X(H_{J_1^*k})\\
                  &= J_1^*\left( \mathcal{L}_X\Pi _1(df,\, dg)\right) - \theta _X(H_{J_1^*k})\\
                  &= -J_1^*\bigl(\phi _1(\xi_f,\, \xi_g,\, X)\bigr) + J_1^*(Xk)\\
                  &= -J_1^*\bigl(\phi _1(\xi_f,\, \xi_g,\, X) + Xk\bigr),  
\end{align*}
where $\xi_f$ means the Hamiltonian vector field of $f\in C^\infty(P_1)$ induced from $\Pi_1$.  
Since $X$ is represented by $X={J_1}_*(H_{J_1^*h})$, we have 
\begin{align*}
Xk &= (dk)\bigl({J_1}_*(H_{J_1^*h})\bigr) = d(J_1^*k)\bigl(H_{J_1^*h}\bigr)\\
   &= -d(J_1^*h)\bigl(H_{J_1^*k}\bigr) = -H_{J_1^*k}(J_1^*h)\\
   & = -J_1^*\phi _1(H_{J_1^*f},\, H_{J_1^*g},\, H_{J_1^*h})\\
   &= -\phi _1(\xi_f,\, \xi_g,\, X)
\end{align*}
Consequently, $d\theta_X=0$ is proved. \qquad\qquad\qquad\qquad\qquad\qquad\qquad\qquad\qquad\qquad\qquad\qquad $\Box$
\subsection{Morita equivalence for twisted symplectic groupoids}

\begin{dfn}\label{twsymplectic}
A twisted symplectic groupoid is a Lie groupoid $\varGamma\rightrightarrows P$ equipped with a non-degenerate 2-form 
$\omega \in \varGamma(\wedge^2T^*\varGamma)$ and a closed 3-form $\phi\in \varGamma(\wedge^3T^*P)$ such that 
 \begin{enumerate}[{\quad \rm (1)}]
  \item $\dim\varGamma = 2\dim P;$
  \item $d\omega \,=\, {\bf t}^*\phi \,-\,{\bf s}^*\phi;$ 
  \item The 2-form $\omega \oplus \omega \oplus (-\omega)$ vanishes on the graph of the 
      groupoid multiplication, 
 \end{enumerate}
\end{dfn}
where ${\bf s}$ and ${\bf t}$ are the source map and target map of $\varGamma$, respectively. 
Before discussing Morita theory for twisted symplectic groupoids, we recall Morita equivalence of Lie groupoids. 
An action of Lie groupoid $G$ on $M$ is said to be 
principal with regard to a smooth map $J:M\rightarrow X$ if $J$ is a surjective submersion and if $G$ acts 
freely and transitively on each $J$-fiber. Assume that Lie groupoids $G\rightrightarrows G_0$ and $H\rightrightarrows H_0$ 
act on a smooth manifold $M$ from the left and right, respectively. We call $M$ a $(G,\,H)$-bibundle 
if the left $G$-action and right $H$-action commute. 
A $(G,\,H)$-bibundle is called a left principal when the left $G$-action is principal with regard to the moment map 
for the right $H$-action. Similarly, it is called a right principal when it is principal with regard to the moment map 
for the left $G$-action. And, a $(G,\,H)$-bibundle is said to be biprincipal when it is both left principal and 
right principal. 
Two Lie groupoids $G$ and $H$ are said to be Morita equivalent if there exists a biprincipal $(G,\,H)$-bibundle. 

\begin{dfn}
Let $(\varGamma\rightrightarrows \varGamma_0,\,\omega_\varGamma,\,\phi_0)$ be a twisted symplectic groupoid. 
A left $\varGamma$-space $M\overset{J}{\rightarrow}\varGamma_0$ is called a twisted symplectic left $\varGamma$-module 
if there exists a non-degenerate $2$-form $\omega_M$ on $M$ such that 
 \begin{enumerate}[{\quad \rm(1)}]
  \item $d\omega_M \,= \,J^*\phi_\varGamma;$
  \item $m^*\omega_M \,=\, {\rm pr}_M^*\omega_M \,+\, {\rm pr}_\varGamma^*\omega_\varGamma,$
 \end{enumerate}
where ${\rm pr}_M$ and ${\rm pr}_\varGamma$ denote the natural projections and 
$m:\varGamma \times_{({\bf s},J)} M\to M$ means the action map. 
\end{dfn}

Twisted symplectic right $\varGamma$-modules are defined in the obvious analogous way. 
We now obtain a category from twisted symplectic $\varGamma$-modules. The category has 
twisted $\varGamma$-modules as objects, and the maps which commute with momentum maps and intertwine 
the $\varGamma$-action as morphisms. We call this category the representation category 
of $\varGamma$ and denote it by ${\rm Rep}(\varGamma)$.

\begin{dfn}\label{Hamiltonianbimodule}
Let $(G\rightrightarrows G_0,\,\omega_G,\,\phi_G)$ and $(H\rightrightarrows H_0,\,\omega_H,\,\phi_H)$ be 
twisted symplectic groupoids. A $(G,\,H)$-bibundle $G_0\overset{\rho}{\leftarrow}X\overset{\sigma}{\rightarrow}H_0$ is 
called a non-degenerate Hamiltonian $(G,\,H)$-bimodule if $X\overset{\rho\times\sigma}{\rightarrow}G_0\times H_0$ is a 
twisted symplectic left $G\times \overline{H}$-module, where the action is given by $(g,\,h)\cdot x=gxh^{-1}$ for 
any $g\in G,,h\in H$ and $x\in X$ such that ${\bf s}(g)=\rho(x)$ and ${\bf s}(h)=\sigma(x)$. 
\end{dfn}

Obviously, the notion of non-degenerate Hamiltonian bimodule is a special case of 
Hamiltonian bimodules defined for quasi-symplectic groupoids in \cite{xu3}. Quasi-symplectic groupoids generalize 
twisted symplectic groupoids. Accordingly, the notion of Morita equivalence for twisted symplectic groupoids is 
given as a part of Morita equivalence of quasi-symplectic groupoids.

\begin{dfn}
Two twisted symplectic groupoids $(G\rightrightarrows G_0,\,\omega_G,\,\phi_G)$ and 
$(H\rightrightarrows H_0,\,\omega_H,\,\phi_H)$ are said to be {\rm (}strong{\rm )} Morita equivalent if they are 
Morita equivalent as Lie groupoids, and the $(G,\,H)$-bibundle is also a non-degenerate Hamiltonian bimodule. 
\end{dfn}

\begin{prop}\label{proposition}
{\rm (\cite{yuji2})}~Assume that $G$ and $H$ are twisted symplectic groupoids that are Morita equivalent 
as Lie groupoids with a Morita bimodule $S$ and a non-degenerate 2-form $\omega_S\in \varGamma(\wedge^2T^*S)$. 
Then, the following conditions are equivalent:
 \begin{enumerate}[{\rm \quad (1)}]
  \item The $G$-action and the $H$-action are both twisted symplectic actions{\rm ;}
  \item $\omega_G\oplus (-\omega_H)\oplus\omega_S \oplus(-\omega_S)$ vanishes on 
        the graph of the $(G\times \overline{H})$-action given in {\rm Definition \ref{Hamiltonianbimodule}}.
 \end{enumerate}
\end{prop}

Due to Proposition \ref{proposition}, twisted symplectic groupoids $G$ and $H$ are Morita equivalence if and only if 
there exists biprincipal $(G,\,H)$-bimodule $G_0\overset{\rho}{\leftarrow}S\overset{\sigma}{\rightarrow}H_0$ equipped 
with a non-degenerate 2-form $\omega_S$ on $S$ such that $d\omega_S=\rho^*\phi_G-\sigma^*\phi_H$ and 
the actions of $G$ and $H$ are twisted symplectic. 
As for symplectic groupoids, the notion of Morita equivalence of twisted symplectic groupoids is closely related to 
Morita equivalence of twisted Poisson manifolds. The following theorem is a generalization of Theorem 3.2 in \cite{xu2}. 

\begin{thm}{\rm (\cite{yuji2})}\label{mainthm}
Let $(P_i,\,\Pi_i,\,\phi_i)~(i=1,2)$ be integrable twisted Poisson manifolds. $P_1$ and $P_2$ are Morita equivalent 
if and only if their associated twisted symplectic groupoid $\mathcal{G}(P_i)\rightrightarrows P_i~(i=1,2)$ are 
Morita equivalent. 
\end{thm}

In order to show Theorem \ref{mainthm}, we need the following result. 
This lemma can be proved in a way similar to Lemma 4.8 
in Bursztyn, H. and Crainic, M. \cite{bursztyn}. 

\begin{lem}\label{lemma}
Let $P$ be a integrable twisted Poisson manifold and $\mathcal{G}(P)\rightrightarrows P$ the associated twisted symplectic 
groupoid. If $S\overset{J}{\to}P$ is a twisted symplectic left $\mathcal{G}(P)$-module, then the momentum map 
$J:S\to P$ is a complete t.s.realization. Conversely, if $J:S\to P$ is a complete t.s.realization, then $S$ is a 
twisted symplectic left $\mathcal{G}(P)$-module. 
\end{lem}

\noindent 
({\it Proof of} Theorem \ref{mainthm})~Assume that $P_1$ and $P_2$ are Morita equivalent 
with an equivalence bimodule $P_1\overset{J_1}{\leftarrow}S\overset{J_2}{\rightarrow}P_2$. 
Let $\gamma : I\to T^*P_1$ be a cotangent path and $y_0\in S$. 
Then, it is known that the horizontal curve $u:I\to S$ over $\gamma$ is given as the solution of 
\[\frac{d}{dt}u(t) = \Pi_1^{\sharp}\bigl(a(t)\bigr),\quad u(0)=y_0\]
and, moreover, depends on the cotangent path $\gamma$~(see \cite{integrability}).
Therefore, for any $g:=[\gamma]\in \mathcal{G}(P_1)$ such that ${\bf s}(g)=x_0$ and ${\bf t}(g)=x_1$, given 
$y_0\in J_1^{-1}(x_0)$, one can define a well-defined action of $\mathcal{G}(P_1)$ on $S$ by 
$g\cdot y_0:= u(1)\in J_1^{-1}(x_1)$. From Lemma \ref{lemma}, this is a twisted symplectic action. 
It is also verified that 
$S\overset{J_2}{\rightarrow}P_2$ has a twisted symplectic action of $\mathcal{G}(P_2)$ in a similar way. 

By assumption and Proposition \ref{eq2}, we have $[H_{J_1^*f},\, H_{J_2^*g}]=0$. This implies that these two actions commute 
(see Theorem 3.2 in \cite{xu2}). 
Therefore, $S$ is a biprincipal $\bigl(\mathcal{G}(P_1),\,\mathcal{G}(P_2)\bigr)$-bibundle. 
By using Proposition \ref{proposition}, it is shown that $\mathcal{G}(P_1)$ and $\mathcal{G}(P_2)$ are Morita equivalent. 

Conversely, we assume that $\mathcal{G}(P_i)\rightrightarrows P_i~(i=1,2)$ are Morita equivalent with a non-degenerate 
Hamiltonian bimodule $P_1\overset{\rho}{\leftarrow}(X,\,\omega_X)\overset{\sigma}{\rightarrow}P_2$. 
Given any $x\in X$, let $g(t)\in \mathcal{G}(P_1)$ and $h(t)\in \mathcal{G}(P_2)$ be smooth curves such that 
${\bf s}\bigl(g(t)\bigr)=\rho(x),\,\varepsilon\bigl(\rho(x)\bigr)=g(0)$ and 
${\bf t}\bigl(h(t)\bigr)=\sigma(x),\, \varepsilon\bigl(\sigma(x)\bigr)=h(0)$, respectively. 
Using $g(t)$ and $h(t)$, we define the new curves $\delta_1(t)$ and $\delta_2(t)$ taking values in the graph 
of the left $\mathcal{G}(P_1)$-action by 
\[
 \delta_1(t) = \bigl(\varepsilon\bigl(\rho(x)\bigr),\,a(t),\,a(t)\bigr),\quad 
 \delta_2(t) = \bigl(g(t),\,x,\,b(t)\bigr), 
\]
where $a(t)=xh(t)$ and $b(t)=g(t)x$. Since $X$ is biprincipal, 
it is verified that $(d\rho)_x(\dot{a}(0))=\boldsymbol{0}$ and $(d\sigma)_x(\dot{b}(0))=\boldsymbol{0}$. 
Accordingly, the $\rho$-fiber and the $\sigma$-fiber coincide with the $\mathcal{G}(P_2)$-orbit and 
the $\mathcal{G}(P_1)$-orbit, respectively. 
If $\omega_G$ is the twisted symplectic form on $\mathcal{G}(P_1)$, then 
\begin{align*}
 0 =&\, -\bigl(\omega_G\oplus\omega_X\oplus(-\omega_X)\bigr)\bigl(\dot{\delta}_1(0),\,\dot{\delta}_2(0)\bigr)\\
   =&\, -\omega_G(\boldsymbol{0},\,\dot{g}(0)) - \omega_X(\dot{a}(0),\,\boldsymbol{0}) + \omega_X(\dot{a}(0),\,\dot{b}(0))\\
   =&\, \omega_X(\dot{a}(0),\,\dot{b}(0))\\
   =&\, i_{\dot{a}(0)}\omega_X, 
\end{align*}
that is, $i_{\dot{a}(0)}\omega_X$ is the element of the annihilation of $\ker(d\sigma)_x$. 
This implies that $\dot{a}(0)$ is represented by the Hamiltonian vector field of some function $F\in C^\infty(P_2)$: 
$\dot{a}(0)=H_{\sigma^*F}$. Similarly, we have $\dot{b}(0)=H_{\rho^*G}~(G\in C^\infty(P_1))$. 
From this observation, it is proved that $\bigl(\ker (d\rho)_x\bigr)^\perp=\ker (d\sigma)_x$ and 
$\bigl(\ker (d\sigma)_x\bigr)^\perp=\ker (d\rho)_x$ for any $x\in X$. 

By assumption, we can show that each the $\rho$-fiber and the $\sigma$-fiber is diffeomorphic to ${\bf s}$-fiber of 
$\mathcal{G}(P_2)$ and ${\bf t}$-fiber of $\mathcal{G}(P_1)$, respectively. Hence, both the $\rho$-fiber and 
the $\sigma$-fiber are connected and simply-connected. Moreover, from Lemma \ref{lemma}, 
$\rho:X\to P_1$ and $\sigma:X\to P_2$ are 
t.s.realization and anti-t.s.realization, respectively. 
Consequently, $P_1$ and $P_2$ are Morita equivalent. 
\qquad\qquad$\Box$
\subsection{Representation categories}

As discussed in \cite{bw}, 
a complete symplectic realization of an integrable Poisson manifold 
plays a role of a left module of algebra. Morita equivalence of Poisson manifolds 
implies equivalent categories of modules, as algebraic Morita equivalence does (\cite{xu2}). 
The result similar to this holds for twisted Poisson manifolds. We define the categories of modules of 
integrable twisted Poisson manifolds as follows: 
\begin{dfn}
Let $P$ be an integrable twisted Poisson manifold. 
The representation category ${\rm Rep}(P)$ of P 
has complete twisted symplectic realizations $S\overset{J}{\to}P$ as objects, and
twisted Poisson maps between twisted symplectic realizations commuting with the realization maps as morphisms. 
\end{dfn}

As already mentioned, for a twisted symplectic groupoid $\varGamma\rightrightarrows\varGamma_0$, 
we can obtain the representation category ${\rm Rep}(\varGamma)$ of $\varGamma$. 
The representation categories are invariants under Morita equivalent.  
That is, the following theorem holds:
\begin{thm}
{\rm (\cite{xu3},\cite{yuji1})}
If two twisted symplectic groupoids $G$ and $H$ are Morita equivalent, then their representation 
categories ${\rm Rep}(G)$ and ${\rm Rep}(H)$ are equivalent. 
\end{thm}

By the above theorem and Theorem \ref{mainthm}, we obtain the following result immediately. 
\begin{thm}{\rm (\cite{yuji1})}
Suppose that $P_1$ and $P_2$ are integrable twisted Poisson manifolds. If $P_1$ and $P_2$ are Morita equivalent, 
then ${\rm Rep}(P_1)$ and ${\rm Rep}(P_2)$ are equivalent. 
\end{thm}

\subsection{Gauge equivalence}
Let $\phi$ be a closed 3-form on a smooth manifold $M$. A $\phi$-twisted Dirac structure on $M$ is a subbundle 
$L_M\subset \mathbb{T}M:=TM\oplus T^*M$ which is maximal isotropic with respect to the symmetric paring 
$\langle\cdot,\,\cdot \rangle$ and whose the set of sections $\varGamma(L_M)$ is closed under 
the bracket $\llbracket\cdot,\,\cdot\rrbracket$, 
where $\langle\cdot,\,\cdot \rangle$ and $\llbracket\cdot,\,\cdot\rrbracket$ are defined as follows:
 \begin{enumerate}[\quad(1)]
  \item $\langle\cdot,\,\cdot \rangle : \varGamma(\mathbb{T}M)\times \varGamma(\mathbb{T}M)\to C^\infty(M),\quad 
   \langle(X,\,\xi),\,(Y,\,\eta)\rangle := \eta (X) + \xi (Y);$
  \item $\llbracket\cdot,\,\cdot\rrbracket : \varGamma(\mathbb{T}M)\times \varGamma(\mathbb{T}M)\to 
     \varGamma(\mathbb{T}M),\quad \llbracket (X,\,\xi),\,(Y,\,\eta)\rrbracket := \bigl(\,[X,\,Y],\,\mathcal{L}_X\eta 
        - i_Yd\xi + i_Xi_Y\phi\,\bigr).$
 \end{enumerate}
For a full discussion of Dirac structures, we refer to \cite{bursztyn}, \cite{bcwz} and 
Bursztyn, H. and Radko, O. \cite{br}. 
\begin{ex}{\rm 
{\rm (}Twisted Poisson manifolds{\rm )}~Let $\phi$ be a closed 3-form on a smooth manifold $P$,and 
$\Pi$ a bivector on $P$. 
The graph $L_\Pi$ of $\Pi^\sharp:T^*P\to TP$ is a $\phi$-twisted Dirac structure if and only if $\Pi$ and $\phi$ satisfy 
the formula (\ref{eq1})~, that is, $P$ is a $\phi$-twisted Poisson manifold. 
}
\end{ex}
\begin{ex}{\rm 
{\rm (}Twisted symplectic manifolds{\rm )}~If $\omega$ is a non-degenerate 2-form on $S$, then the graph $L_\omega$ 
of $\omega^\flat:TS\to T^*S$ is a $\psi$-twisted Dirac structure if and only if $\omega$ is a $\psi$-twisted 
symplectic form. 
}
\end{ex}

Let $(M,\,L_M,\,\phi_M)$ and $(N,\,L_N,\,\phi_N)$ be twisted Dirac manifolds. A smooth map $J:M\to N$ is said to be a 
forward Dirac map if 
\[(L_N)_{J(x)} \,=\, \left\{\bigl((dJ)_xV,\,\alpha\bigr)\mid V\in T_xM,\,\alpha\in T_{J(x)}N,\,
                                             \left(V,\,(dJ)_x^*\alpha\right)\in (L_M)_x\right\}\]
for all $x\in M$. 
We write $\bigl(\mathfrak{F}(J)\bigr)(L_M)$ for the right-hand side in the above formula. 
As verified easily, if $L_M$ and $L_N$ are associated with twisted Poisson structures, then a forward Dirac map is 
equivalent to a twisted Poisson map. 

\vspace*{0.5cm}
Now we recall  gauge transformations on ($\phi$-)twisted Dirac structures. Let $L_M$ be a $\phi$-twisted Dirac structure 
on $M$ and $B\in \varGamma(\wedge^2T^*M)$. We set 
\[\tau_B(L_M):=\left\{\bigl(X,\,\xi + B^\flat(X)\bigr)\bigm| (X,\,\xi)\in L_M\right\}.\]
The subbundle $\tau_B(L_M)$ defines a $(\phi - dB)$-twisted Dirac structure on $M$. The operation 
$L_M\mapsto \tau_B(L_M)$ is called a gauge transformation of $L_M$ associated with $B$. Especially, if 
$L_\Pi$ is a $\phi$-Dirac structure associated with a twisted Poisson manifold $P$, 
the gauge transformation associated with $B$ is given by 
\[L_\Pi\longmapsto \tau_B(L_\Pi)=\left\{\bigl(\Pi^\sharp(\alpha),\,\alpha + B^\flat(\Pi^\sharp(\alpha))\bigr)\bigm| 
\alpha\in T^*P\right\}.\]
As discussed in \cite{severa}, $\tau_B(L_\Pi)$ may fail to be induced from a twisted Poisson bivector. 
$\tau_B(L_\Pi)$ is associated with a $(\phi-dB)$-twisted Poisson manifold if and only if 
$1+B^\flat\Pi^\sharp:T^*P\to T^*P$ is invertible. Two twisted Poisson manifold $(P,\,\Pi,\,\phi)$ and $(P,\,\Pi',\,\phi')$ 
are said to be gauge equivalent if there exists a 2-form $B$ on $P$ such that 
\[\Pi' = \Pi\circ (1+B^\flat\Pi^\sharp)^{-1}\quad \text{and}\quad \phi-\phi'=dB.\]
For a twisted Poisson bivector which is gauge equivalent to $P$ with respect to $B$, we write $\tau_B(\Pi)$. 
\begin{thm}
Let $(P,\,\Pi,\,\phi)$ be an integrable twisted Poisson manifold and $(\mathcal{G}(P),\,\omega)$ the associated twisted 
symplectic groupoid. 
Then, for any 2-form $B\in \varGamma(\wedge^2T^*P)$ such that 
$(1+B^\flat\Pi^\sharp)$ is invertible, 
$(P,\,\Pi,\,\phi)$ and $(P,\,\tau_B(\Pi),\,\phi -dB)$ are Morita equivalent with a equivalence bimodule 
$(\mathcal{G}(P),\,\hat{\omega})$, where $\hat{\omega}:=\omega - {\bf s}^*B$. 
\end{thm}
({\it Proof}\,)~For any $x\in \mathcal{G}(P)$, we set $V=T_x\bigl(\mathcal{G}(P)\bigr),\,W=T_{{\bf s}(x)}P$. We define 
$H_1$ and $H_2$ by
\begin{align*}
  &H_1:= \tau_{{\bf s}^*B}(L_\omega) = \{\,(v,\,i_v(\omega + {\bf s}^*B))\mid v\in V\,\}, \\
  &H_2:= \bigl(\mathfrak{F}({\bf s})\bigr)(L_\omega) = \{\,({\bf s}_*v,\,\eta) 
                                    \mid v\in V,\,\eta\in W^*,\,\eta\circ{\bf s}=i_v\omega\,\}.
 \end{align*} 
Then, 
\begin{align*}
  \bigl(\mathfrak{F}({\bf s})\bigr)(H_1) &= \{\,(d{\bf s}(v),\,\eta)
                          \mid v\in V,\,\eta\in W^*,\,\eta\circ{\bf s}=i_v(\omega + {\bf s}^*B)\,\}\\
  \tau_B(H_2) &= \{\,(d{\bf s}(v),\,\xi + i_{{\bf s}_*v}B)\mid v\in V,\,\xi\in W^*,\,{\bf s}^*\xi=i_v\omega\,\}\\
              &= \{\,(d{\bf s}(v),\,\eta)\mid v\in V,\,\eta\in W^*,\,{\bf s}^*(\eta - i_{{\bf s}_*v}B)=i_v\omega\,\}. 
 \end{align*}
Since $i_v({\bf s}^*B)={\bf s}^*\bigl(i_{{\bf s}_*v}B\bigr)$, we have ${\bf s}^*\eta=i_v(\omega + {\bf s}^*B)$. 
Therefore, $\bigl(\mathfrak{F}({\bf s})\bigr)(H_1) = \tau_B(H_2)$. 
From the fact that ${\bf s}:\mathcal{G}(P)\to (P,\,-\Pi)$ is a twisted Poisson map~(\cite{cx}), it follows that 
$\bigl(\mathfrak{F}({\bf s})\bigr)(H_1) = \tau_B(-\Pi)$, that is, 
${\bf s}:(\mathcal{G}(P),\,\widehat{\omega})\to (P,\,\tau_B(\Pi))$ is an anti-t.s.realization. 
It also can be shown that ${\bf t}:(\mathcal{G}(P),\,\widehat{\omega})\to (P,\,\Pi)$ is a t.s.realization in a similar way. 
Moreover, using $(\ker d{\bf s})^\omega = \ker d{\bf t}$, we have
\begin{align*}
   (\ker d{\bf s})^{\widehat{\omega}} &= \{\,v\in V\mid \widehat{\omega}(v,\,w)=0\;(\forall w\in \ker d{\bf s})\,\}\\
                  &= \{\,v\in V\mid \omega(v,\,w)=0\;(\forall w\in \ker d{\bf s})\,\} = (\ker d{\bf s})^\omega\\
                  &= \ker d{\bf t}.
 \end{align*} 
Similarly, we can prove $(\ker d{\bf t})^{\widehat{\omega}}=\ker d{\bf s}$. 

In what follows, we prove that $\bf t$ and $\bf s$ are complete. Let $\hat{H}_\bullet$ and $\hat{X}_\bullet$ denote 
the Hamiltonian vector field with regard to $\hat{\omega}$ and $\omega'$, respectively. 
Then, from assumption we have ${\bf s}^*B(H_{{\bf t}^*f}) = \boldsymbol{0}~(f\in C^\infty(P))$. 
Therefore, $\hat{\omega}({H}_{{\bf t}^*f},\,\cdot) = \omega({H}_{{\bf t}^*f},\,\cdot)=d({\bf t}^*f)(\cdot)$. 
This implies that $\hat{H}_{{\bf t}^*f} = H_{\bf t}^*f$. From the completeness of 
${\bf t}:(\mathcal{G}(P),\,\omega)\to P$, we can conclude that ${\bf t}:(\mathcal{G}(P),\,\hat{\omega})\to P$ is complete. 
The completeness of ${\bf s}$ can be proved similarly (see \cite{br}), 
\qquad\qquad\qquad\qquad\qquad\qquad\qquad\qquad\qquad $\Box$

\section{Weak Morita equivalence}

Let $A_i\to M_i~(i=1,2)$ be Lie algebroids. Assume that $A_1$ and $A_2$ act on $X$ 
from the left and right, respectively. 
If the actions $\varrho_1,\,\varrho_2$ commute i.e., $[\varrho_1(\xi),\,\varrho_2(\eta)]=0$ 
for any $\xi\in \varGamma(A_1),\,\eta\in \varGamma(A_2)$, 
and the moment maps are surjective submersions, then we call $X$ an $(A_1,\,A_2)$-algebroid bimodule. 

\begin{dfn}
Let $(P_i,\,\Pi_i,\,\phi_i)~(i=1,2)$ be twisted Poisson manifolds. $P_1$ and $P_2$ are said to be (weak) Morita 
equivalent if there exists a $(T^*P_1,\,T^*P_2)$-algebroid bimodule $P_1\overset{J_1}{\leftarrow}M\overset{J_2}{\to}P_2$ 
which satisfies 
 \begin{enumerate}[\rm \quad(1)]
  \item Each $J_i$-fiber~$(i=1,2)$ is connected and simply-connected{\rm ;}
  \item For any $x\in M$, 
       \[T_x\bigl(J_1^{-1}(J_1(x))\bigr) = \left\{\,\varrho_2(\eta)_x\bigm| \eta\in \varGamma(T^*P_2)\,\right\}\quad\text{and} 
         \quad T_x\bigl(J_2^{-1}(J_2(x))\bigr) = \left\{\,\varrho_1(\xi)_x\bigm| \xi\in \varGamma(T^*P_1)\,\right\}, \]
  where $\varrho_1:\varGamma(T^*P_1)\to \varGamma(TM)$ and $\varrho_2:\varGamma(T^*P_2)\to \varGamma(TM)$ mean 
  the Lie algebroid actions of $T^*P_1$ and $T^*P_2$, respectively. 
 \end{enumerate}
\end{dfn}

\begin{thm}
Weak Morita equivalence is an equivalence relation for (twisted) Poisson manifolds. 
\end{thm}
({\it Proof}\,)~First, we verify the reflectivity. If $P$ is a twisted Poisson manifold, its cotangent bundle 
$T^*P\overset{\pi}{\to}P$ is a $(T^*P,\,T^*P)$-algebroid bimodule 
under the left action 
$\varrho_L:\varGamma(T^*P)\to \varGamma\bigl(T(T^*P)\bigr),\,\alpha\mapsto \Pi_C^\sharp(\pi^*\alpha)$ and 
the right action $\varrho_R:\varGamma(T^*P)\to \varGamma\bigl(T(T^*P)\bigr),\, \alpha\mapsto -\Pi_C(\pi^*\alpha)$, 
where $\Pi_C$ means a Poisson bivector induced from a canonical symplectic structure on $T^*P$.
From $\pi_*\Pi_C^\sharp(\pi^*\alpha)=0$, we have $T_u\bigl(\pi^{-1}(\pi(u))\bigr) = 
\rho_L\bigl(\varGamma(T^*P)\bigr)_u = \rho_R\bigl(\varGamma(T^*P)\bigr)_u~(\forall u\in T^*P)$. 
Therefore, $P$ is weak Morita equivalent to itself. 

As for the symmetry, we suppose that $P_1$ is weak Morita equivalent to $P_2$ with an algebroid bimodule 
$P_1\overset{J_1}{\leftarrow}M\overset{J_2}{\to}P_2$. Then, a smooth manifold 
$P_2\overset{J_2}{\leftarrow}M\overset{J_1}{\to}P_1$ with the reversed actions $\varrho'_L:=-\varrho_R,\,
\varrho'_R:=-\varrho_L$ is a $(T^*P_2,\,T^*P_1)$-algebroid bimodule. It follows that $P_2$ is weak Morita equivalent 
to $P_1$. 

The transitivity will be shown in what follows. Suppose that $P_1\overset{J_1}{\leftarrow}M\overset{J_2}{\to}P_2$ is a 
$(T^*P_1,\,T^*P_2)$-algebroid bimodule 
and $P_2\overset{J'_2}{\leftarrow}N\overset{J_3}{\to}P_3$ is a $(T^*P_2,\,T^*P_3)$-algebroid bimodule. 
We define the left and right actions on 
the fiber product $L:=M\times^{J_2,J_2'}_{P_2} N$ by 
\[\tilde{\varrho}_1:\varGamma(T^*P_1)\to \varGamma(TL),\, \alpha\mapsto\bigl(\varrho_1(\alpha),\,\boldsymbol{0}\bigr)
\quad\text{and}\quad
\tilde{\varrho}_3:\varGamma(T^*P_3)\to \varGamma(TL),\,\beta\mapsto \bigl(\boldsymbol{0},\, \varrho_3(\beta)\bigr),\]
respectively. 
Then, $L$ is a $(T^*P_1,\,T^*P_3)$-algebroid bimodule with the moment maps 
$P_1\overset{\rho}{\leftarrow}L\overset{\sigma}{\to}P_3$, where $\rho(m,\,n):=J_1(m)$ and $\sigma(m.\,n):=J_3(n)$. 
From assumption, we have 
\[T_{(m,n)}\left(\sigma^{-1}\bigl(\sigma(m,\,n)\bigr)\right) 
      = T_{(m.n)}\left(J_2^{-1}\bigl(J'_2(n)\bigr)\times \{n\}\right)
              =\tilde{\varrho}_1\bigl(\varGamma(T^*P_1)\bigr)_{(m,n)}.\]
Similarly, $T_{(m,n)}\left(\rho^{-1}\bigl(\rho(m,\,n)\bigr)\right) 
                 = \tilde{\varrho}_3\bigl(\varGamma(T^*P_3)\bigr)_{(m,n)}$. Obviously, each fiber is connected and 
simply-connected. Hence, $P_1$ and $P_3$ is weak Morita equivalent. \qquad\qquad\qquad\qquad\qquad\qquad\qquad\quad $\Box$

\begin{prop}
Strong Morita equivalence implies weak Morita equivalence. 
\end{prop}
({\it Proof}\,)~Assume that integrable twisted Poisson manifolds $P_1$ and $P_2$ are strong Morita equivalent with 
an equivalence bimodule $P_1\overset{J_1}{\leftarrow}(S,\,\omega_S)\overset{J_2}{\rightarrow}P_2$. 
As discussed in Section 2, the moment maps $J_1$ and $J_2$ induce Lie algebroid actions 
$\varrho_1(\alpha):= \Pi_S^\sharp(J_1^*\alpha)$ and $\varrho_2(\beta):= \Pi_S(J_2^*\beta)$, respectively, 
where $\Pi_S$ means the bivector field induced from $\omega_S$. 
Using (\ref{eq2}), we have $[\varrho_1(\alpha),\,\varrho_2(\beta)]=0$ for any $\alpha\in 
\varGamma(T^*P_1),\,\beta\in \varGamma(T^*P_2)$. 
This implies that $S$ is a $(T^*P_1,\,T^*P_2)$-algebroid bimodule. 
From assumption, it follows that, for any $x\in S$,
\[T_x\left(J_1^{-1}\bigl(J_1(x)\bigr)\right) = \ker(dJ_1)_x = \bigl(\ker(dJ_2)_x\bigr)^\perp = 
\Pi_S^\sharp\bigl(\ker(dJ_2)_x^\circ\bigr) = \varrho_2\bigl(\varGamma(T^*P_2)\bigr)_x,\]
where $\ker(dJ_2)_x^\circ$ denotes the annihilator of $\ker(dJ_2)_x$. 
Analogously, we have 
$T_x\left(J_2^{-1}\bigl(J_2(x)\bigr)\right) = \varrho_1\bigl(\varGamma(T^*P_1)\bigr)_x$. 
Therefore, $P_1$ and $P_2$ are weak Morita equivalent. \qquad\qquad\qquad\qquad\qquad\qquad $\Box$

Weak Morita equivalence induces one-to-one correspondence between twisted symplectic leaves. 
The following theorem can be shown in a way similar to Theorem 11.1.9 in Ortega, J. and Ratiu, T. \cite{ortega}. 

\begin{thm}\label{leaf}
Suppose that $P_1$ and $P_2$ are weak Morita equivalent with a algebroid bimodule 
$P_1\overset{J_1}{\leftarrow}M\overset{J_2}{\rightarrow}P_2$. Let $M/\mathcal{D}$ be the leaf space of the distribution 
$\mathcal{D}$ defined by $\mathcal{D}_m:=\ker (dJ_1)_m + \ker (dJ_2)_m\,(\forall m\in M)$ and 
$\mathcal{L}(P_i)~(i=1,2)$ the spaces of twisted symplectic leaves of $P_i$, respectively. 
Then, 
 \begin{enumerate}[{\quad\rm (1)}]
  \item The distribution $\mathcal{D}=\ker (dJ_1) + \ker (dJ_2)$ is integrable.
  \item $M/\mathcal{D}\to \mathcal{L}(P_i)~(i=1,2)$ are bijections. In particular, the map 
  $\mathcal{L}(P_1)\to \mathcal{L}(P_2),\,L\mapsto J_2\bigl(J_1^{-1}(L)\bigr)$ is the bijective correspondence between 
  the leaves of $P_1$ and the leaves of $P_2$.
 \end{enumerate}
\end{thm}

\noindent ({\it Proof}\,)~ (1)~
From assumption, $\ker (dJ_1)$ and $\ker (dJ_2)$ can be considered as the distributions  
\[V_1=\{\,\varrho_2(\beta)\,|\,\beta\in \varGamma(T^*P_2)\,\}\quad \text{and}\quad 
V_2=\{\,\varrho_1(\alpha)\,|\,\alpha\in \varGamma(T^*P_1)\,\},\] 
respectively. 
The distribution $\mathcal{D}$ is spanned by $V=V_1\cup V_2$. 
Let $\theta_t$ and $\eta_t$ be the flows of $\varrho_1(\alpha)$ and $\varrho_2(\beta)$, respectively. 
We will show that 
\[(d\theta_t)_m\bigl(\varrho_2(\beta)_m\bigr)\in \mathcal{D}\bigl(\theta_t(m)\bigr)\quad\text{and}\quad
 (d\eta_t)_m\bigl(\varrho_1(\alpha)_m\bigr)\in \mathcal{D}\bigl(\eta_t(m)\bigr).\]
Since $\theta_t$ is a diffeomorphism, we can define the pull-back of $X\in \varGamma(TM)$ by 
$\theta_t^*X:=(d\theta_{-t})\circ X\circ \theta_t$. 
Then, we have the following formula~(see \cite{ortega}):
\[\frac{d}{dt}\theta_t^*\bigl(\varrho_2(\beta)\bigr) = \theta_t^*[\varrho_1(\alpha),\,\varrho_2(\beta)].\]
By assumption that the two Lie algebroid actions commute, the right-hand side in the above formula is equal 
to $\boldsymbol{0}$. Therefore, $\theta_t^*\varrho_2(\beta) = \theta_0^*\varrho_2(\beta)=\varrho_2(\beta)$.
It follows that 
\[(d\theta_t)_m\bigl(\varrho_2(\beta)_m\bigr) = 
    (d\theta_t)_m\circ (d\theta_{-t})_{\theta_t(m)}\bigl(\varrho_2(\beta)\bigr)_{\theta_t(m)}
     = \bigl(\varrho_2(\beta)\bigr)_{\theta_t(m)}\in \mathcal{D}\bigl(\theta_t(m)\bigr).\]
We can show that $(d\eta_t)_m\bigl(\varrho_1(\alpha)_m\bigr)\in \mathcal{D}\bigl(\eta_t(m)\bigr)$ in similar way. 

\vspace*{0.5cm}
\noindent (2)~
We will denote by $\Phi_V,\,\Phi_{V_1}$ and $\Phi_{V_2}$ the pseudogroups of local transformations generated 
by the flows of elements in $V,\, V_1$ and $V_2$, respectively. For full discussion of pseudogroups, 
we refer to \cite{ortega}. 
Let $N\subset M$ be the integrable manifold of $\mathcal{D}$ containing a given point $m\in M$. 
We note that $L$ coincides with the $\Phi_V$-orbit of $m$: 
\[N=\Phi_V\cdot m = \{\,\varphi(m)\mid \varphi\in \Phi_V\,\}.\] 
Since the two Lie algebroid actions commute, $N$ can be written as 
\[N = \Phi_V\cdot m = \Phi_{V_1}\bigl(\Phi_{V_2}\cdot m\bigr).\]
From assumption, the $J_i$-fibers $(i=1,2)$ are preserved by the elements in $\Phi_{V_i}$, respectively. Accordingly, 
\[J_1(N) = J_1\bigl(\Phi_{V_1}(\Phi_{V_2}\cdot m)\bigr) = J_1(\Phi_{V_2}\cdot m).\]
Any element $\theta\in \Phi_{V_2}$ can be represented as $\theta=\theta^1_{t_1}\circ \cdots \circ\theta^n_{t_n}$, where 
$\theta^j_{t_j}~(j=1,\cdots,n)$ mean the flows of a vector fields $\varrho_1(dJ_1^*f_j),\,f_j\in C^\infty(P_1)$. 
Accordingly, 
\[J_1\bigl(\theta(m)\bigr) = J_1\bigl((\theta^1_{t_1}\circ \cdots \circ\theta^n_{t_n})(m)\bigr)
 = \bigl(\xi^1_{t_1}\circ \cdots \circ\xi^n_{t_n}\bigr)\bigl(J_1(m)\bigr),\]
where $\xi^j_{t_j}~(j=1,\cdots,n)$ are the flows of Hamiltonian vector fields $H_{f_j}$. 
$\xi^1_{t_1}\circ \cdots \circ\xi^n_{t_n}$ is the element of the pseudogroup $\Phi_H$ of local transformations 
generated by the flows of Hamiltonian vector fields on $P_1$. 
Moreover, the twisted symplectic leaf of $P_1$ is the maximal integral manifold of the distribution spanned by 
the Hamiltonian vector fields on $P_1$. 
Therefore, we have $J_1(N)=J_1(\Phi_{V_2}\cdot m)=L_{J_1(m)}$, 
where $L_{J_1(m)}$ means the leaf of $P_1$ containing $J_1(m)$.
Consequently, we can define the map $\Psi:M/\mathcal{D}\rightarrow \mathcal{L}(P_1)$ by
\[\Psi: M/\mathcal{D}\longrightarrow \mathcal{L}(P_1),\quad N=\Phi_V\cdot m\longmapsto J_1(N)=\Phi_H\cdot J_1(m).\]
To show the bijectivity of $\Psi$, 
we will prove that the map $\Psi'$ defined by $J_1(N)\mapsto J_1^{-1}\bigl(J_1(N)\bigr)$ is an inverse of $\Psi$. 
From $J_1(N)=J_1(\Phi_{V_2}\cdot m)$, it follows that  
\[J_1^{-1}\bigl(J_1(N)\bigr) = \bigcup_{\theta\in \Phi_{V_2}}J_1^{-1}\bigl(J_1(\theta(m))\bigr).\]
Here, since the elements in $\Phi_{V_1}$ preserve each $J_1$-fiber, we have 
$J_1^{-1}\bigl(J_1(\theta(m))\bigr) = \Phi_{V_1}\cdot\theta(m)$
for any $\theta\in \Phi_{V_2}$. Therefore, 
\[J_1^{-1}\bigl(J_1(N)\bigr) = \bigcup_{\theta\in \Phi_{V_2}}\Phi_{V_1}\cdot\theta(m)
=\Phi_{V_1}\cdot(\Phi_{V_2}\cdot m)=\Phi_V\cdot m=N.\]
This leads us to the conclusion that $\Psi$ is bijective. Similarly, we can construct the map 
$M/\mathcal{D}\to \mathcal{L}(P_2)$ and show that this map is bijective. 
\qquad\qquad\qquad\qquad\qquad\qquad\qquad\qquad\qquad\qquad\qquad $\Box$

\begin{minipage}[t]{10cm}
Yuji HIROTA\\
Quantum Bio-Informatics Center\\ 
Research Institute for Science and Technology\\
Tokyo University of Science\\
2641 Yamazaki, Noda City, Chiba\\ 
278-8510 Japan;\\
email:yhirota@rs.noda.tus.ac.jp
\end{minipage}

\end{document}